\documentclass{amsart}
\usepackage{amssymb}
\usepackage{graphicx}

\newtheorem{thm}{Theorem}
\newtheorem{cor}{Corollary}
\newtheorem{lem}{Lemma}

\newtheorem{rem}{Remark}

\newcommand{\es}{{\mathcal S}}

\newcommand{\D}{{\mathbb D}}

\def\be{\begin{equation}}
\def\ee{\end{equation}}

\newcommand{\bee}{\begin{enumerate}}
\newcommand{\eee}{\end{enumerate}}

\newcommand{\blem}{\begin{lem}}
\newcommand{\elem}{\end{lem}}
\newcommand{\bthm}{\begin{thm}}
\newcommand{\ethm}{\end{thm}}
\newcommand{\bcor}{\begin{cor}}
\newcommand{\ecor}{\end{cor}}
\newcommand{\beg}{\begin{example}}
\newcommand{\eeg}{\end{example}}
\newcommand{\begs}{\begin{examples}}
\newcommand{\eegs}{\end{examples}}
\newcommand{\bdefe}{\begin{defin}}
\newcommand{\edefe}{\end{defin}}
\newcommand{\bprob}{\begin{prob}}
\newcommand{\eprob}{\end{prob}}
\newcommand{\bei}{\begin{itemize}}
\newcommand{\eei}{\end{itemize}}

\newcommand{\bcon}{\begin{conj}}
\newcommand{\econ}{\end{conj}}
\newcommand{\bcons}{\begin{conjs}}
\newcommand{\econs}{\end{conjs}}
\newcommand{\bprop}{\begin{propo}}
\newcommand{\eprop}{\end{propo}}
\newcommand{\br}{\begin{rem}}
\newcommand{\er}{\end{rem}}
\newcommand{\brs}{\begin{rems}}
\newcommand{\ers}{\end{rems}}
\newcommand{\bo}{\begin{obser}}
\newcommand{\eo}{\end{obser}}
\newcommand{\bos}{\begin{obsers}}
\newcommand{\eos}{\end{obsers}}
\newcommand{\bpf}{\begin{pf}}
\newcommand{\epf}{\end{pf}}
\newcommand{\ba}{\begin{array}}
\newcommand{\ea}{\end{array}}
\newcommand{\beq}{\begin{eqnarray}}
\newcommand{\beqq}{\begin{eqnarray*}}
\newcommand{\eeq}{\end{eqnarray}}
\newcommand{\eeqq}{\end{eqnarray*}}

\begin{document}
\bibliographystyle{amsplain}

\title[On the difference of coefficients of univalent functions]{On the difference of coefficients of univalent functions}

\author[M. Obradovi\'{c}]{Milutin Obradovi\'{c}}
\address{Department of Mathematics,
Faculty of Civil Engineering, University of Belgrade,
Bulevar Kralja Aleksandra 73, 11000, Belgrade, Serbia}
\email{obrad@grf.bg.ac.rs}

\author[D. K. Thomas]{Derek K. Thomas}
\address{Department of Mathematics, Swansea University, Bay Campus, Swansea, SA1 8EN, United Kingdom}
\email{d.k.thomas@swansea.ac.uk}

\author[N. Tuneski]{Nikola Tuneski}
\address{Department of Mathematics and Informatics, Faculty of Mechanical Engineering, Ss. Cyril and Methodius
University in Skopje, Karpo\v{s} II b.b., 1000 Skopje, Republic of North Macedonia.}
\email{nikola.tuneski@mf.edu.mk}





\begin{abstract}
For $f\in \mathcal{S}$, the class of normalized functions, analytic and univalent in the unit disk $\mathbb{D}$ and given by $f(z)=z+\sum_{n=2}^{\infty} a_n z^n$ for $z\in \mathbb{D}$, we give an upper bound for the coefficient difference $|a_4|-|a_3|$  when $f\in \mathcal{S}$. This  provides an improved bound in the case $n=3$  of Grispan's 1976 general bound $||a_{n+1}|-|a_n||\le 3.61\dots .$ Other coefficients bounds, and bounds for the second and third Hankel determinants when $f\in \mathcal{S}$ are found when either $a_2=0,$ or $a_3=0$.
\end{abstract}

\maketitle

\medskip

\section{Introduction. preliminaries and definitions}

\medskip

Let $\mathcal{A}$ be the class of functions $f$ which are analytic  in the open unit disc $\D=\{z:|z|<1\}$ of the form
\be\label{e1}
f(z)=z+a_2z^2+a_3z^3+\cdots,
\ee
and let $\mathcal{S}$ be the subclass of $\mathcal{A}$ consisting of functions that are univalent in $\D$.

\medskip

Although the famous Bieberbach conjecture $|a_n|\le n$ for $n\ge2,$ was proved by de Branges  in 1985 \cite{Bra85}, a great many other problems concerning the coefficients $a_n$ remain open. The main aim of this paper (Section 3), is by use of  the Grunsky inequalities, to find an upper for the difference of coefficients $|a_4|-|a_3|$ for $f\in \mathcal{
S}$, which improves the well-known general bound of Grispan $||a_{n+1}|-|a_n||\le 3.61\dots$  \cite{Gri76}, when $n=3$.
We also obtain information concerning the initial coefficients of $f(z)$, and of the second and third Hankel determinants  when either $a_2=0$, or $a_3=0$.
\bigskip

For $f\in \mathcal{S}$, the Grunsky coefficients $\omega_{p,q}$ as defined in N. A. Lebedev \cite{Lebedev} are given by
\[  \log\frac{f(t)-f(z)}{t-z}=\sum_{p,q=0}^{\infty}\omega_{p,q}t^{p}z^{q},  \]

\medskip

\noindent where  $\omega_{p,q}=\omega_{q,p}$, and satisfy the so-called Grunsky inequalities \cite{duren,Lebedev}

\be\label{eq 4}
\sum_{q=1}^{\infty}q \left|\sum_{p=1}^{\infty}\omega_{p,q}x_{p}\right|^{2}\leq \sum_{p=1}^{\infty}\frac{|x_{p}|^{2}}{p},
\ee
where $x_{p}$ are arbitrary complex numbers such that last series converges.

\medskip

Further, it is well-known that if $f$ given by \eqref{e1}
belongs to $\mathcal{S}$, then also
\be\label{eq 5-mo-4}
f_{2}(z)=\sqrt{f(z^{2})}=z +c_{3}z^3+c_{5}z^{5}+\cdots
\ee
belongs to $\mathcal{S}$. Thus for the function $f_{2}$ we have the appropriate Grunsky
coefficients of the form $\omega_{2p-1,2q-1}$, and  inequalities \eqref{eq 4} take the form
\be\label{eq 6-mo-5}
\sum_{q=1}^{\infty}(2q-1) \left|\sum_{p=1}^{\infty}\omega_{2p-1,2q-1}x_{2p-1}\right|^{2}\leq \sum_{p=1}^{\infty}\frac{|x_{2p-1}|^{2}}{2p-1}.
\ee

\medskip

(Note that in this paper, we omit the upper index (2) in  $\omega_{2p-1,2q-1}^{(2)}$ in Lebedev's notation).

\medskip

The following similar inequality follows from the relation (15) on page 57 in \cite{Lebedev}.
\be\label{eq 6-mo-6}
 \left|\sum_{p=1}^{\infty}\sum_{q=1}^{\infty}\omega_{2p-1,2q-1}x_{2p-1}x_{2q-1}\right|\leq \sum_{p=1}^{\infty}\frac{|x_{2p-1}|^{2}}{2p-1}.
\ee

\medskip

Thus for example, from  \eqref{eq 6-mo-5} and \eqref{eq 6-mo-6} when $x_{2p-1}=0$ and $p=3,4,\ldots$, we obtain
\begin{equation}\label{e7}
|\omega_{11} x_1 +\omega_{31} x_3 |^2 +3|\omega_{13} x_1 +\omega_{33} x_3 |^2 + 5|\omega_{15} x_1 +\omega_{35} x_3 |^2 \le |x_1|^2+\frac{|x_3|^2}{3}
\end{equation}
and
\begin{equation}\label{e8}
|\omega_{11} x_1^2 +2\omega_{13}x_1 x_3 +\omega_{33}x_3^2 |\le |x_1|^2+\frac{|x_3|^2}{3},
\end{equation}
respectively.

\medskip

It was also shown in \cite[p.57]{Lebedev}, that if $f\in\mathcal{S}$ is given by \eqref{e1}, then the coefficients $a_{2}$, $ a_{3}$, $ a_{4}$ and $a_5$ can be expressed in terms of the Grunsky coefficients  $\omega_{2p-1,2q-1}$ of the function $f_{2}$ given by
\eqref{eq 5-mo-4} as follows.
\be\label{e9}
\begin{split}
a_{2}&=2\omega _{11},\\
a_{3}&=2\omega_{13}+3\omega_{11}^{2}, \\
a_{4}&=2\omega_{33}+8\omega_{11}\omega_{13}+\frac{10}{3}\omega_{11}^{3},\\
a_{5}&=2\omega_{35}+8\omega_{11}\omega_{33}+5\omega_{13}^{2}+18\omega_{11}^2\omega_{13}+\frac73\omega_{11}^4,\\
0&= 3\omega_{15}-3\omega_{11}\omega_{13}+\omega_{11}^3-3\omega_{33}.
\end{split}
\ee
In this paper we will use these expressions to obtain information concerning the coefficients $a_2$, $a_3$, $a_4$, and $a_5$ when $f\in\mathcal{S}$.
\medskip

In recent years a great deal of attention has been given to finding upper bounds for the modulus of the second and third Hankel determinants $H_2(2)$ and $H_3(1)$, defined as follows who's elements are the coefficients of $f\in\mathcal{S}$ (see e.g. \cite{DTV-book}).
 \medskip

For $f\in\mathcal{S}$
\[  H_2(2) = a_2a_4-a_3^2  \]
and
\begin{equation}\label{e11}
H_3(1) = a_3(a_2a_4-a_3^2) - a_4(a_4-a_2a_3)+a_5(a_3-a_2^2).
\end{equation}
\medskip

Almost all results have concentrated on finding bounds for $|H_2(2)|$ and $|H_3(1)|$ for subclasses of $\mathcal{S}$, and only recently has a significant bound been found for the whole class $\mathcal{S}$ \cite{OT} for $|H_2(2)|$ and $|H_3(1)|$. However finding exact sharp bounds remains an open problem.
\medskip

We begin by using the Grunsky inequalities in (\ref{eq 6-mo-6}) to obtain bounds for the modulus of some initial coefficients and $|H_2(2)|$ and $|H_3(1)|$ when  $f\in \mathcal{S}$ provided either $a_2,$ or $a_3=0.$

\medskip
\section{Coefficient bounds and Hankel determinants}

\medskip

Obtaining sharp bounds for the modulus of the coefficients for odd functions in $\mathcal{S}$ has long been been an open problem. If $f_2$, given by (\ref{eq 5-mo-4}) is an odd function in $\mathcal{S}$, then  the only known sharp bounds for $|c_{2n-1}|$ for $n\ge 2$ are $|c_3|\le 1$, and $|c_5|\le 1/2+e^{-2/3}=1.013\dots$. In general the best bound to date is $|c_{2n-1}|\le1.14 $ for $n\ge2$, (see e.g.\cite{duren}).
\medskip

In our first theorem, we give  bounds for $|a_3|$, $|a_4|$ and $|a_5|$ when $f\in \mathcal{S}$   assuming only that only $a_2=0,$ thus providing bounds for a wider class of functions than the odd functions in $\mathcal{S}$. We also give bounds for $|H_2(2)|$ and $|H_3(1)|$  in this case.
\begin{thm}\label{th-1}
Let $f\in\es$ and be given by \eqref{e1} with $a_2=0$. Then
\begin{itemize}
\item[($i$)] $|a_3|\le1$,
\medskip
\item[($ii$)] $|a_4|\le\frac23=0.666\ldots$,
\medskip
\item[($iii$)]  $|a_5|\le \sqrt{\frac{19}{15}} = 1.67666\ldots$,
\medskip
\item[($iv$)] $|H_2(2)|\le1$,
\medskip
\item[($v$)] $|H_3(1)|\le\frac{21}{20} = 1.05$.
\end{itemize}
\end{thm}
\medskip
\begin{proof}$\mbox{}$
\begin{itemize}
\item[($i$)] The classical inequality $|a_3-a_2^2|\le1$ for $f$ in $\es$ when $a_2=0$, gives  $|a_3|\le1$, which from (\ref{e9}) gives
\be\label{e12}
|\omega_{13}|\le\frac12.
\ee
\medskip
\item[($ii$)] Next choose  $x_1=0$ and $x_3=1$ in \eqref{e8}, which gives
\be\label{e13}
|\omega_{33}| \le\frac13.
\ee
Also, since $\omega_{11}=0$ ($\Leftrightarrow a_2=0$), then from \eqref{e9} and \eqref{e13} we obtain
\[
|a_4|=2|\omega_{33}|\le\frac23=0.666\ldots.
\]
\medskip
\item[($iii$)]  Again since  $\omega_{11}=0$,  from \eqref{e9} we obtain
\be\label{e14}
|a_5|=|2\omega_{35}+5\omega_{13}^2|.
\ee

From \eqref{e7} with $x_1=0$ and $x_3=1$ we have ($\omega_{11}=0$)
\[  |\omega_{13}|^2 + 3|\omega_{33}|^2 + 5|\omega_{35}|^2 \le \frac13 \]
and from here
\begin{equation}\label{n1}
  |\omega_{35}| \le \frac{1}{\sqrt{15}} \sqrt{1-3|\omega_{13}|^2}.
\end{equation}
From \eqref{e14} and \eqref{n1} we have
\[ |a_{5}| \le 2|\omega_{35}|+5|\omega_{13}|^2 \le \frac{1}{\sqrt{15}} \sqrt{1-3|\omega_{13}|^2}+5|\omega_{13}|^2 \le \frac{503}{300}=1.67666\ldots.\]

\medskip

\item[($iv$)] Since we are assuming  $a_2=0$, $(i)$ shows that   $|H_2(2)| \le1$ is trivial.
\medskip

\item[($v$)] When $\omega_{11}=0$, from the last relation in \eqref{e9} we have $\omega_{33}=\omega_{15},$ and from \eqref{e11},
\be\label{e16}
\begin{split}
|H_3(1)
&= |2\omega_{13}^3+4\omega_{13}\omega_{35}-4\omega_{33}^2|\le 2|\omega_{13}|^3 + 4 + \underbrace{|\omega_{13}\omega_{35} -\omega_{15}^2|}_{E_1}.
\end{split}
\ee

Now choose $x_1=-\omega_{15}$, and $x_3=\omega_{13}$, and since $\omega_{33}=\omega_{15}$, from \eqref{e7} we obtain
\[ |\omega_{13}|^4+5E_1^2\le |\omega_{15}|^2+\frac{|\omega_{13}|^2}{3} \le \frac15-\frac35|\omega_{13}|^2  +\frac13|\omega_{13}|^2,\]
(since by \eqref{e7} $3|\omega_{13}|^2+5|\omega_{15}|^2\le1$ for $x_1=1$, $x_3=0$ and $\omega_{11}=0$),
which implies $5E_1^2\le \frac15-\frac{4}{15}|\omega_{13}|^2-|\omega_{13}|^4$, i.e., $E_1\le \frac15$.

\medskip

Finally from \eqref{e12} and \eqref{e16}, it follows that
\[ |H_3(1)| \le 2\cdot \frac18+4\cdot \frac15 = \frac{21}{20}=1.05. \]
\end{itemize}
This completes the proof of Theorem \ref{th-1}.
\end{proof}

\medskip

We next prove a similar result, this time assuming that $a_3=0$.
\begin{thm}\label{th-2}
Let $f\in\es$ and be  given by \eqref{e1}, with $a_3=0$. Then
\medskip
\begin{itemize}
\item[($i$)] $|a_2|\le1$,
\medskip
\item[($ii$)] $|a_4|\le \frac{\sqrt{37}+13}{12}=1.59023\ldots$,
\medskip
\item[($iii$)]  $|a_5|\le \frac14 \sqrt{\frac{757}{15}}+ \frac{85}{64} = 3.10412\ldots$,
\medskip
\item[($iv$)] $|H_2(2)|\le\frac{13+\sqrt{37}}{12}=1.59023\ldots$,
\medskip
\item[($v$)] $|H_3(1)|\le \frac{24+\sqrt{645}}{30}  = 1.64656\ldots$.
\end{itemize}
\end{thm}
\medskip

\begin{proof}$\mbox{}$
\begin{itemize}
\item[($i$)] Since $|a_3-a_2^2|\le1$ and $a_3=0$, then $|a_2^2|\le1$, i.e., $|a_2|\le1$. Also, since by \eqref{e9}, $a_3=2\omega_{13}+3\omega_{11}^2=0$, it follows that
\be\label{e17}
\omega_{13}=-\frac32 \omega_{11}^2 \quad \left(\Leftrightarrow\,\, \omega_{11}^2=-\frac23\omega_{13}\right).
\ee
Because $|a_2|=|2\omega_{11}|\le1$, we have
\be\label{e18}
|\omega_{11}|\le \frac12 \quad\mbox{and}\quad |\omega_{13}|\le \frac38 \,\,(\mbox{by } \eqref{e17}.
\ee
\medskip
\item[($ii$)] By using \eqref{e9} and \eqref{e17}, we obtain
\be\label{e19}
\begin{split}
|a_4|
&= \left|2\omega_{33}+8\omega_{11} \left(-\frac32\omega_{11}^2\right)+\frac{10}{3}\omega_{11}^3\right|\\
&= \left|2\omega_{33}-\frac{26}{3}\omega_{11}^3\right| \\
&\le 2|\omega_{33}|+\frac{26}{3}|\omega_{11}| ^3.
\end{split}
\ee
From \eqref{e7}, using $x_1=0$ and $x_3=1$, we have
\[  |\omega_{13}|^2 + 3|\omega_{33}|^2 \le \frac13, \]
which implies (with $\omega_{13}=-\frac32 \omega_{11}^2$, see \eqref{e17})
\be\label{e20}
|\omega_{33}|\le\sqrt{\frac19-\frac34 |\omega_{11}|^4}.
\ee
Combining \eqref{e19} and \eqref{e20} we obtain
\be\label{e21}
|a_4|\le 2 \sqrt{\frac19 -\frac34|\omega_{11}|^4} + \frac{26}{3}|\omega_{11}|^3 =: \varphi(|\omega_{11}|),
\ee
where $\varphi(t) =  2 \sqrt{\frac19 -\frac34 t^4} + \frac{26}{3}t^3 $, $0\le t =|\omega_{11}|\le\frac12$ (by \eqref{e18}). Since $\varphi$ is increasing function on  $[0,1/2]$,
\[ \varphi(t)\le \varphi(1/2)=\frac{\sqrt{37}+13}{12}, \]
which, together   with \eqref{e21}, gives the desired result.
\medskip
\item[($iii$)] From the last relation in \eqref{e9}, using \eqref{e17}  we have $\omega_{33}=\omega_{15}+\frac{11}{6}\omega_{11}^3$, which with the expression for $a_5$ in \eqref{e9}, gives
\be\label{e22}
\begin{split}
|a_5|
&= |2\omega_{35} + 8\omega_{11}\omega_{15} +5\omega_{13}^2 -10\omega_{11}^4|\\
&\le 2\underbrace{| \omega_{35} + 4\omega_{11}\omega_{15}|}_{C_1^\ast} + \underbrace{5|\omega_{13}|^2 +10|\omega_{11}|^4}_{C_2^\ast}.
\end{split}
\ee
Once again, using \eqref{e7} choosing $x_1=4\omega_{11}$, $x_3=1$ and $\omega_{13}=-\frac32 \omega_{11}^2$, we have
\[(C_1^\ast)^2 = |4\omega_{11}\omega_{15} +\omega_{35} |^2 \le -\frac54 |\omega_{11}|^4 +\frac{16}{5}|\omega_{11}|^2+\frac{1}{15} \le \frac{757}{64\cdot 15}, \]
since $|\omega_{11}|\le\frac12$. Thus
\[ C_1^\ast \le \frac18\sqrt{\frac{757}{15}}.\]

\medskip

Next, since $\omega_{13}=-\frac32 \omega_{11}^2$ and $|\omega_{11}|\le\frac12$, we have
\[
\begin{split}
C_2^\ast
&= 5\cdot\frac94\cdot |\omega_{11}|^4+10|\omega_{11}|^4 = \frac{85}{4}|\omega_{11}|^4 \le \frac{85}{4}\cdot \frac{1}{16} = \frac{85}{64},
\end{split}
\]
since $|\omega_{11}|\le\frac12$.

\medskip

Finally from  \eqref{e22} we have
\[  |a_5| \le \frac14 \sqrt{\frac{757}{15}}+ \frac{85}{64} = 3.10412\ldots.\]
\medskip
\item[($iv$)] By using \eqref{e11}, \eqref{e9} and \eqref{e17}, we have
\be\label{e24}
\begin{split}
H_2(2)
&= 4\omega_{11}\omega_{33} + 4\omega_{11}^2\omega_{13} - 4\omega_{13}^2 - \frac73\omega_{11}^4\\
&= 4\omega_{11}\omega_{33}-\frac{52}{3}\omega_{11}^4
\end{split}
\ee
and from here
\be\label{n2}
\begin{split}
|H_2(2)| \le 4|\omega_{11}||\omega_{33}|+\frac{52}{3}|\omega_{11}|^4.
\end{split}
\ee

From \eqref{e20} and \eqref{n2} we have
\[ |H_2(2)| \le 4|\omega_{11}|\sqrt{\frac19-\frac34|\omega_{11}|^4} +\frac{52}{3}|\omega_{11}|^4 =: \varphi_1(|\omega_{11}|,\]
where
\[ \varphi_1(t) = 4t\sqrt{\frac19-\frac34t^4}+\frac{52}{3}t^4, \]
with $0\le t=|\omega_{11}|\le\frac12$. Finally, it can be checked that  $\varphi_1$ is an increasing function on the interval $(0,1/2)$, and so
\[ |H_2(2)| \le \varphi_1(1/2) = \frac{13+\sqrt{37}}{12}=1.59023\ldots. \]

\medskip
\item[($v$)]
By using the last relation from \eqref{e9} with $\omega_{13}=-\frac32\omega_{11}^2$, it follows that $\omega_{33}=\omega_{15}+\frac{11}{6}\omega_{11}^3$, and so using  \eqref{e11}, after some calculations we obtain
\[ H_3(1) = -12\omega_{11}^2 \left(\omega_{11}\omega_{15} +\frac23\omega_{35}\right) - 4\omega_{15}^2 -30\omega_{11}^6,\]
which gives
\be\label{e25}
|H_3(1)| \le \underbrace{12 |\omega_{11}|^2 \left|\omega_{11}\omega_{15} +\frac23 \omega_{35}\right|}_{D_1} + \underbrace{4|\omega_{15}|^2  + 30|\omega_{11}|^6}_{D_2}.
\ee
Now choose $x_1=\omega_{11}$ and $x_3=\frac23$ in \eqref{e7}, then (since $\omega_{13}=-\frac32\omega_{11}^2$),
\[  \left|\omega_{11}\omega_{15} +\frac23\omega_{35}\right| \le \sqrt{\frac15\left( |\omega_{11}|^2+\frac{4}{27} \right)},\]
and so
\be\label{e26}
\begin{split}
D_1 &\le 12 |\omega_{11}|^2 \sqrt{\frac15\left(|\omega_{11}|^2+\frac{4}{27}\right)} \\
&\le 12 \cdot\frac14 \sqrt{\frac15\left(\frac14+\frac{4}{27}\right)} = \sqrt{\frac{43}{60}}=\frac{\sqrt{645}}{30}=0.84656\dots,
\end{split}
\ee
since $|\omega_{11}|\le\frac12$.

\medskip

Also, as in the proof of ($iii$), we have
\[ 5|\omega_{15}|^2\le 1-|\omega_{11}|^2 -3|\omega_{13}|^2 = 1-|\omega_{11}|^2-\frac{27}{4}|\omega_{11}|^4,  \]
where we have once again used $\omega_{13}=-\frac32\omega_{11}^2$.
Now
\[ D_2 \le \frac45 - \frac{4}{5}|\omega_{11}|^2 - \frac{27}{5}|\omega_{11}|^4 + 30|\omega_{11}|^6  =: \varphi_2(|\omega_{11}|^2), \]
where
\[  \varphi_2(t)=\frac15\left(4 -4t - 27t^2 +150t^3\right),  \]
and $0\le t = |\omega_{11}|^2 \le \frac14$. Since $\varphi_2$ attains its maximum at $t_0 = 0$,
\be\label{e29}
D_2 \le \varphi_2(0) = \frac45.
\ee
Finally, by using \eqref{e25}, \eqref{e26} and \eqref{e29} we obtain
\[ |H_3(1)| \le D_1+D_2 \le  \frac{24+\sqrt{645}}{30}  = 1.64656\ldots. \]
\end{itemize}
\end{proof}

\medskip

\section{Coefficient  differences for $f\in \mathcal{S}$}
\medskip

A long standing problem in the theory of univalent functions is to find  sharp upper and lower bounds for $|a_{n+1}|-|a_n|$, when $f\in\mathcal{S}$. Since  the Keobe function has coefficients $a_n=n$, it is natural to conjecture  that $||a_{n+1}|-|a_n||\le1$. As early as 1933, this was shown to be false even when $n=2,$ when Fekete and Szeg\"o \cite{FekSze33} obtained the sharp bounds
 $$ -1 \leq |a_3| - |a_2| \leq \frac{3}{4} + e^{-\lambda_0}(2e^{-\lambda_0}-1) = 1.029\ldots, $$
where $\lambda_0$ is the unique value of $\lambda$ in $0 < \lambda <1$, satisfying the equation $4\lambda = e^{\lambda}$.

\medskip

Hayman \cite{Hay63} showed that if $f \in {\mathcal S}$, then $| |a_{n+1}| - |a_n| | \leq C$, where $C$ is an absolute constant. The exact value of $C$ is unknown,  the best estimate to date
being $C=3.61\ldots$ \cite{Gri76}, which because of the sharp estimate above when $n=2$, cannot be reduced to $1$.

\medskip

We now use the methods of this paper to obtain a better upper bound in the case $n=3$.

\medskip

\begin{thm}
Let $f\in\es$ and be given by \eqref{e1}. Then
\[|a_4|-|a_3| \le 2.1033299\ldots.\]
\end{thm}

\medskip

\begin{proof}
 By using \eqref{e9} we have
\[|a_4|-|a_3| \le |a_4|-|\omega_{11}||a_3| \le  |a_4-\omega_{11} a_3| = 2\Big|\underbrace{\omega_{33} +3\omega_{11}\omega_{33} + \frac16 \omega_{11}^3}_{B}\Big|. \]
From \eqref{e8} with $x_1=\frac{1}{\sqrt6}\omega_{11}$ and $x_3=1$, we obtain
\[
\begin{split}
 & \left| \omega_{33} +\frac{2}{\sqrt6}\omega_{11}\omega_{13} +\frac16\omega_{11}^3 \right| \le \frac16 |\omega_{11}|^2 +\frac13  \\
 \Rightarrow \quad & \left|B+\left( \frac{2}{\sqrt6}-3 \right) \omega_{11}\omega_{13} \right| \le \frac16 |\omega_{11}|^2 +\frac13  \\
 \Rightarrow \quad & |B| \le \left(3- \frac{\sqrt6}{3} \right) |\omega_{11}||\omega_{13}| + \frac16 |\omega_{11}|^2 +\frac13  \\
 \Rightarrow \quad & |B| \le \left(3- \frac{\sqrt6}{3} \right) |\omega_{11}| \cdot \frac{1}{\sqrt3}\sqrt{1-|\omega_{11}|^2} + \frac16 |\omega_{11}|^2 +\frac13  \\
 \Rightarrow \quad & |B| \le \frac13\left[ (3\sqrt3-\sqrt2) |\omega_{11}| \sqrt{1-|\omega_{11}|^2} + \frac12 |\omega_{11}|^2 + 1 \right] =: \varphi(|\omega_{11}|),
\end{split}
\]
where $\varphi(t)= \frac13\left[ (3\sqrt3-\sqrt2) t \sqrt{1-t^2} + \frac12 t^2 + 1 \right]$ for $0\le t\le1$, and where we have used that $|\omega_{13}| \le \frac{1}{\sqrt3}\sqrt{1-|\omega_{11}|^2}$. Since the function $\varphi$ attains its maximum at
\[t_0 = \sqrt{\frac12+\frac16\sqrt{\frac{1}{379}(39+8\sqrt6)}} = 0.75202\ldots,\]
and since $\varphi(t_0) = \frac{1}{12}\left(5+\sqrt{117-24\sqrt6}\right)$, it follows that
\[ |a_4|-|a_3| \le 2\varphi(t_0) = 2.10495\ldots. \]
\end{proof}

\end{document}